\documentclass{article}
\usepackage{amssymb}
\usepackage{amsmath}
\usepackage{bbm}
\usepackage{geometry}
\usepackage{authblk}
\usepackage{amsthm}

\theoremstyle{plain}
\newtheorem{theorem}{Theorem}[section]
\newtheorem{proposition}{Proposition}[section]
\newtheorem{lemma}{Lemma}[section]

\theoremstyle{definition}
\newtheorem{definition}{Definition}[section]
\newtheorem{claim}{Claim}[section]

\theoremstyle{remark}
\newtheorem{remark}{Remark}[section]

\newcommand{\conv}{\textup{conv}}

\title{Equality conditions for the fractional superadditive volume inequalities}
\author{Mark Meyer\footnote{Department of Mathematics, Wichita State University, Wichita, KS, USA, Email: memeyer@shockers.wichita.edu}}
\date{July 14, 2023}
\begin{document}
\maketitle
\begin{abstract}
While studying set function properties of Lebesgue measure, F. Barthe and M. Madiman proved that Lebesgue measure is fractionally superadditive on compact sets in $\mathbb{R}^n$. In doing this they proved a fractional generalization of the Brunn-Minkowski-Lyusternik (BML) inequality in dimension $n=1$. In this paper we will prove the equality conditions for the fractional superadditive volume inequalites for any dimension. The non-trivial equality conditions are as follows. In the one-dimensional case we will show that for a fractional partition $(\mathcal{G},\beta)$ and nonempty sets $A_1,\dots,A_m\subseteq\mathbb{R}$, equality holds iff for each $S\in\mathcal{G}$, the set $\sum_{i\in S}A_i$ is an interval. In the case of dimension $n\geq2$ we will show that equality can hold if and only if the set $\sum_{i=1}^{m}A_i$ has measure $0$.
\end{abstract}
\textbf{Classification Code.} 52A40\\
\textbf{Data Availability Statement.}Data sharing not applicable to this article as no datasets were generated or analysed during the current study.

\vspace{0.8cm}

\section{Introduction}

\noindent The Brunn-Minkowski inequality is a well known result of convex geometry. Due to H. Brunn and H. Minkowski, the result states for compact, convex sets $A, B$ in $\mathbb{R}^n$ and $\lambda\in[0,1]$ that 
\begin{equation*}
|\lambda A+(1-\lambda)B|^{\frac{1}{n}}\geq\lambda|A|^{\frac{1}{n}}+(1-\lambda)|B|^{\frac{1}{n}},
\end{equation*}
where $|A|$ represents $n$-dimensional Lebesgue measure of a set $A\subset\mathbb{R}^n$. The result was later extended by L. Lyusternik in 1935 to compact sets, and the resulting inequality is called the Brunn-Minkowski-Lyusternik (BML) inequality.

\begin{theorem}[The Brunn-Minkowski-Lyusternik inequality \cite{Lusternik1}]\label{BML_inequality} Let $A$ and $B$ be non-empty compact sets in $\mathbb{R}^n$. Then 
\begin{equation*}
|A+B|^{\frac{1}{n}}\geq|A|^{\frac{1}{n}}+|B|^{\frac{1}{n}},
\end{equation*}
where equality holds if and only if exactly one of the following conditions hold:
\begin{enumerate}
    \item $|A+B|=0$.
    \item Exactly one of the sets $A,B$ is a point the other has positive measure.
    \item The sets $A$ and $B$ are full dimensional homothetic convex bodies. (If $n=1$, then this condition states that $A$ and $B$ are intervals with positive measure.)
\end{enumerate}
\end{theorem}

Some useful references to the Brunn-Minkowski inequality are \cite{Burago1}, which provides a proof of the equality conditions for the BML inequality and \cite{Gardner1}, which provides a survey about the Brunn-Minkowski inequality.

What follows is a summary of some of the motivation for the work in this paper.

\textbf{1. The fractional Young's inequality and the fractional Brunn-Minkowski inequality:} One result of the BML inequality is that the functional $A\rightarrow|A|^{\frac{1}{n}}$ is superadditive on compact sets in $\mathbb{R}^n$. This means that if $A_1,\dots,A_m$ are non-empty compact sets in $\mathbb{R}^n$, then
\begin{equation*}
\left|\sum_{i=1}^{m}A_i\right|^{\frac{1}{n}}\geq\sum_{i=1}^{m}|A_i|^{\frac{1}{n}}.
\end{equation*}
Bobkov et al. conjectured in \cite{Bobkov1} that there is a fractional generalization of both Young's convolution inequality and the Brunn-Minkowski inequality, and it was shown that the former implies the latter. In particular, it was conjectured that the above functional is fractionally superadditive. That is, if $A_1,\dots,A_m$ are compact sets in $\mathbb{R}^n$, and if $(\mathcal{G},\beta)$ is a fractional partition over $\{1,\dots,m\}$ (see Section \ref{notation_definitions}), then
\begin{equation*}
\left|\sum_{i=1}^{m}A_i\right|^{\frac{1}{n}}\geq\sum_{S\in\mathcal{G}}\beta(S)\left|\sum_{i\in S}A_i\right|^{\frac{1}{n}}.
\end{equation*}
We are interested in the generalization of the Brunn-Minkowski inequality. If the fractional Brunn-Minkowski inequality is true, it would give an extensive generalization of the Brunn-Minkowski inequality when the number of sets being added is larger than 2. Since fractional superadditivity is an important property in the study of set functions (see the appendix in \cite{Barthe1} for a discussion on this), this kind of extension for the Brunn-Minkowski inequality would be interesting. It was immediately verified in \cite{Bobkov1} that the fractional superadditivity conjecture is true when the sets are compact and convex. But, if the sets under consideration are compact, but not necessarily convex, the fractional BML inequality may not be true. In \cite{Fradelizi1} it was shown that the fractional BML inequality is not true in dimension 12 or higher. Later in \cite{Fradelizi3} it was shown that the conjecture fails in dimension $n\geq7$. It has been verified in \cite{Barthe1} that the fractional BML conjecture is true in dimension 1.

\textbf{2. The monotonicity of the volume deficit:} In \cite{Fradelizi1,Fradelizi2} there are some developments related to the fractional superadditive BML inequality. The general idea is to observe the set average of a non-empty compact set $A\subset\mathbb{R}^n$ defined by $A(k)=\frac{1}{k}\sum_{i=1}^{k}A$ and to observe how close $A(k)$ is to $\text{conv}(A)$. The Shapely-Folkmann-Starr theorem \cite{Starr1} states that $A(k)$ converges to $\text{conv}(A)$ in Hausdorff distance. It has also been shown in \cite{Emerson1} that the volume deficit $\Delta(A):=|\text{conv}(A)\backslash A(k)|$ converges to 0 as $k\rightarrow\infty$ (as long as $|A(k)|>0$ for some $k$). A natural question to ask is whether the volume deficit is non-increasing in $k$. In \cite{Fradelizi2} the volume deficit is determined to be non-increasing in dimension 1, but not necessarily for dimension 12 or higher. Later in \cite{Fradelizi3} it was determined that for star shaped sets the volume deficit is monotone in $k$ for dimensions 2 and 3, and also for dimension $n\geq 4$ as long as $k\geq(n-1)(n-2)$. The connection between fractional superadditivity for the functional $A\rightarrow |A|^{\frac{1}{n}}$ and the monotonicity of the volume deficit can be seen by observing the fractional superadditive inequality for the leave-one-out partition. The leave-one-out partition for the integers $\{1,\dots,m\}$ is the partition that consists of the sets $S_j:=\{1,\dots,m\}\backslash\{j\}$, where each set $S_j$ is assigned a weight of $\beta(S_j)=1/(m-1)$. Then the corresponding fractional superadditive inequality conjectured by Bobkov et al. is
\begin{equation*}
    \left|\sum_{i=1}^{m}A_i\right|^{\frac{1}{n}}\geq\frac{1}{m-1}\sum_{j=1}^{m}\left|\sum_{i\in S_j}A_i\right|^{\frac{1}{n}}.
\end{equation*}
If we set $A_i=A$, then this becomes 
\begin{equation*}
\left|\sum_{i=1}^{m}A\right|^{\frac{1}{n}}\geq\frac{1}{m-1}\sum_{i=1}^{m}\left|\sum_{j=1}^{m-1}A\right|^{\frac{1}{n}}=\frac{m}{m-1}\left|\sum_{i=1}^{m-1}A\right|^{\frac{1}{n}}.
\end{equation*}
By rearranging terms and using the n-homogenous property of Lebesgue measure we find that the above inequality is just $|A(m)|\geq|A(m-1)|$ which then implies that the volume deficit is non-increasing. It follows that if the fractional superadditive BML inequality is true, then the volume deficit is monotone in $k$. On the other hand, if the volume deficit is not non-increasing for some set $A\subset\mathbb{R}^n$, a counter example has been produced which shows that the fractional superadditive BML inequality for the leave-one-out partition is not true.

\textbf{3. The Shannon-Stam entropy power inequality:} There is study of the analogy between the functional $f_X\rightarrow N(X)$, where $X$ is a random vector with pdf $f_X$, $N(X)$ is the power entropy of $X$, and the functional $A\rightarrow|A|^{\frac{1}{n}}$, where $A$ is compact in $\mathbb{R}^n$. A result in information theory known as the Shannon-Stam entropy power inequality \cite{Shannon1,Stam1} is analogous to the BML inequality. The Shannon-Stam entropy power inequality shows that the power entropy $N(X)$ of a random variable $X$ has a superadditive property: $N(X+Y)\geq N(X)+N(Y)$ if $X$ and $Y$ are independent random vectors. Later the fractional superadditive property of the power entropy was determined (\cite{Artstein1,Madiman1,Madiman2}).

\textbf{4. Supermodularity properties of volume:} There is an interest in determining certain properties of the set function $v(S)=|\sum_{i\in S}A_i|$, $S\subset [m]$. In addition to fractional superadditivity, another property of interest when studying set functions is supermodularity. The function $v$ is supermodular if for any $S,T\subset[m]$ it is true that
\begin{equation*}
v(S\cup T)+v(S\cap T)\geq v(S)+v(T).
\end{equation*}
There is a connection between these two properties. If the set function $v$ is supermodular, then $v$ is also fractionally superadditive \cite{ollagnier1}. When the sets $A_1,\dots,A_m$ are compact and convex, the function $v$ is supermodular \cite{Fradelizi2}. This was determined by proving the inequality
\begin{equation*}
|A+B+C|+|A|\geq|A+B|+|A+C|,
\end{equation*}
which implies that $v$ is supermodular. When the sets $A_1,\dots,A_m$ are compact and not necessarily convex this inequality was found to be false. As a result the supermodularity of $v$ could not be verified. But a similar inequality was verified for dimension 1. If $A,B,C$ are compact sets in $\mathbb{R}$, then
\begin{equation*}
|A+B+C|+|\text{conv}(A)|\geq|A+B|+|A+C|.
\end{equation*}
In \cite{Fradelizi4} this inequality was verified to be true with certain restrictions when the dimension is greater than 1. 

 While it is known that the fractional BML inequality is not always true it would be interesting to know if at least the volume functional has the fractional superadditive property when the sets are compact. It would also be interesting to see that even though the function $v(S)=|\sum_{i\in S}A_i|$ is not supermodular on compact sets, that $v$ is fractionally superadditive. Recent work in \cite{Barthe1} has shown that Lebesgue measure does have the fractional superadditive property, proving that the fractional BML inequality in \cite{Bobkov1} is true for $n=1$ when the sets are compact.

\begin{theorem}[Barthe and Madiman \cite{Barthe1}]\label{Fractional_Superadditive}
For an integer $m\geq3$ let $(\mathcal{G},\beta)$ be a fractional partition of $[m]$. Then
\begin{equation*}
\left|\sum_{i=1}^{m}A_i\right|\geq\sum_{S\in\mathcal{G}}\beta(S)\left|\sum_{i\in S}A_i\right|
\end{equation*}
holds for any non-empty compact subsets $A_1,\dots,A_m$ of $\mathbb{R}^n$.
\end{theorem}

We emphasize that Theorem \ref{Fractional_Superadditive} establishes fractional superadditivity for the volume functional $|A|$, which establishes the conjectured BML inequality in \cite{Bobkov1} without the power $1/n$. In \cite{Barthe1} it was given as an open question to completely characterize the conditions for equality when $A_1,\dots,A_m$ are compact subsets of $\mathbb{R}$. Barthe and Madiman made initial progress toward resolving this conjecture by proving that equality holds (in dimension $1$) if all sets $S\in\mathcal{G}$ with $\beta(S)>0$ satisfy $\sum_{i\in S}A_i$ is an interval. The objective of this paper is to provide a complete characterization of the equality conditions in any dimension. In the case of dimension $1$ we will prove that the converse statement also holds. We will prove that

\begin{theorem}\label{dimension_1_equality}
    Let $m\geq3$ be an integer, let $(\mathcal{G},\beta)$ be a fractional partition of $[m]$, and let $A_1,\dots,A_m$ be compact sets in $\mathbb{R}$. Then
    \begin{equation*}
        \left|\sum_{i=1}^{m}A_i\right|=\sum_{S\in\mathcal{G}}\beta(S)\left|\sum_{i\in S}A_i\right|
    \end{equation*}
    if and only one of the following conditions holds:
    \begin{enumerate}
        \item $\left|\sum_{i\in[m]}A_i\right|=0$.
        \item $\left|\sum_{i\in[m]}A_i\right|>0$ and the translated fractional partition (see Remark \ref{R2}) is trivial.
        \item $\left|\sum_{i\in[m]}A_i\right|>0$,the translated fractional partition is non-trivial, and for each $S\in\mathcal{G}$ with $\beta(S)>0$ the set $\sum_{i\in S}A_i$ is an interval (where a point is considered as an interval).
    \end{enumerate}
\end{theorem}

It turns out that in the case of dimension $n\geq2$ there is an extra degree of freedom which eliminates the convexity condition in Theorem \ref{dimension_1_equality}. We will show 
\begin{theorem}\label{dimension_n_equality}
     Let $n,m\geq2$ be integers, let $A_1,\dots,A_m$ be nonempty compact sets in $\mathbb{R}^n$, and let $(\mathcal{G},\beta)$ be a fractional partition of $[m]$. Then
    \begin{equation*}
        \left|\sum_{i=1}^{m}A_i\right|=\sum_{S\in\mathcal{G}}\beta(S)\left|\sum_{i\in S}A_i\right|
    \end{equation*}
    if and only if one of the conditions is true:
    \begin{enumerate}
        \item The translated fractional partition is trivial.
        \item The translated fractional partition is non-trivial and $\left|\sum_{i=1}^{m}A_i\right|=0$.
    \end{enumerate}
\end{theorem}

The remainder of this paper will be dedicated to proving the above theorems and a few other interesting facts.

\subsection{Acknowledgements}

I would like to thank my PhD advisor Robert Fraser for introducing me to topics in Convex Geometry and for the assistance he has given me as I worked on this problem. In particular, he suggested it might be useful to look carefully at the proof of Theorem \ref{Fractional_Superadditive}. I am grateful that he has taken the time to read various editions of this article and offered helpful corrections and advice.  I thank the anonymous referees for their suggestions which have helped to significantly improve the quality of this paper.

\section{Notation and definitions}\label{notation_definitions}

For the following let $A$ represent a compact set in $\mathbb{R}$. If $a\in\mathbb{R}$, then define ${}_{a\leq} (A):=[a,\infty)\cap A$ and $(A)_{\leq a}:=(-\infty,a]\cap A$. 
To denote the points $x\in A$ which are surrounded by positive measure in $A$ we use
\begin{equation*}
\text{PM}(A):=\{x\in A: |(x-\varepsilon,x+\varepsilon)\cap A|>0\text{   } \forall \varepsilon>0\}.
\end{equation*}
Since the set $\text{PM}(A)$ may not be closed, the set $\overline{\text{PM}(A)}$ will be used to denote its closure. For $x\in\text{conv}(A)\backslash\overline{\text{PM}(A)}$ define the closest point of $\overline{\text{PM}(A)}$ to the right of $x$ by
\begin{equation*}
x_R:=\inf\{y\in\overline{\text{PM}(A)}:y>x\},
\end{equation*}
and if $\{y\in\overline{\text{PM}(A)}:y>x\}=\varnothing$, then we say $x_R$ does not exist. Similarly define the closest point of $\overline{\text{PM}(A)}$ to the left of $x$ by
\begin{equation*}
x_L:=\sup\{y\in\overline{\text{PM}(A)}:y<x\},
\end{equation*}
and if $\{y\in\overline{\text{PM}(A)}:y<x\}=\varnothing$, then we say $x_L$ does not exist. In the case that $A$ is compact in $\mathbb{R}^n$, then define
\begin{equation*}
    \textup{PM}(A):=\{x\in\mathbb{R}^n:|B(x,\varepsilon)\cap A|>0\text{  }\forall \varepsilon>0\},
\end{equation*}
where $B(x,\varepsilon):=\{x\in \mathbb{R}^n:|x|\leq\varepsilon\}$ is the closed ball with center $x$ and radius $\varepsilon$. 

For the rest of this work, unless stated otherwise, it will be assumed, because of the translation invariance of Lebesgue measure, that $\text{min}(A_i)=0$ and that $\text{max}(A_i)=a_i\geq0$ for any compact set $A_i\subseteq\mathbb{R}$. Now we formally define a fractional partition.

\begin{definition} A \textit{fractional partition} over $[m]:=\{1,\dots,m\}$ is a pair $(\mathcal{G},\beta)$ where $\mathcal{G}$ is a collection of subsets of $[m]$ and $\beta:2^{[m]}\rightarrow[0,1]$ is a function of weights on $\mathcal{G}$ which satisfies: $\forall i\in[m]$, $\sum_{S\in\mathcal{G}:i\in S}\beta(S)=1$.
\end{definition}

To avoid ambiguity three additional requirements will be added to the definition of a fractional partition. The first requirement is that $\beta(S)>0$ for each $S\in\mathcal{G}$. We can make this assumption since removing sets $S\in\mathcal{G}$ with $\beta(S)=0$ has no effect on the fractional superadditive inequalities in Theorem \ref{Fractional_Superadditive}. The second requirement is that if $(\mathcal{G},\beta)$ is not the trivial partition, then $[m]\notin\mathcal{G}$. To justify this requirement assume that $[m]\in\mathcal{G}$ and that $(\mathcal{G},\beta)$ is not the trivial partition so that $\beta([m])<1$. Then the fractional superadditive inequality is
\begin{equation*}
\left|\sum_{i=1}^{m}A_i\right|\geq\beta([m])\left|\sum_{i=1}^{m}A_i\right|+\sum_{S\in\mathcal{G}\backslash[m]}\beta(S)\left|\sum_{i\in S}A_i\right|.
\end{equation*}
This can be simplified to get
\begin{equation*}
\left|\sum_{i=1}^{m}A_i\right|\geq\sum_{S\in\mathcal{G}\backslash[m]}\frac{\beta(S)}{1-\beta([m])}\left|\sum_{i\in S}A_i\right|.
\end{equation*}
We now observe that for each $S\in\mathcal{G}\backslash[m]$
\begin{equation*}
\sum_{S\in\mathcal{G}\backslash[m]:i\in S}\frac{\beta(S)}{1-\beta([m])}=\frac{1-\beta([m])}{1-\beta([m])}=1.
\end{equation*}
This shows that every nontrivial fractional partition over $[m]$ corresponds to a fractional superadditive inequality for which $\beta([m])|\sum_{i=1}^{m}A_i|$ does not belong to the right side of the inequality. It follows that the term $\beta([m])|\sum_{i=1}^{m}A_i|$ on the right side of the inequality has no influence over the equality conditions. The third requirement is that the weights $\beta(S)$ are rational. This assumption allows us to reduce the inequality of Theorem \ref{Fractional_Superadditive} to the inequality $q|\sum_{i=1}^{m}A_i|\geq\sum_{j=1}^{s}|\sum_{i\in S_j}A_i|$, where $q\geq1$ is an integer, the sets $S_j$ are an enumeration of the sets $S\in\mathcal{G}$ possibly with repetition, and each $i\in[m]$ belongs to exactly $q$ of the sets $S_j$. The reason that this assumption can be made is complicated and the justification is given in Section \ref{rational weights}

\section{Initial observations}\label{initial_observations}

The objective of this section is to provide an outline of the proof of Theorem \ref{Fractional_Superadditive}. We will point out what kind of modifications can be made without causing the proof to break down. The observations made in this section are crucial in understanding the proofs of the following results. Recall the assumption that for each $A_1,\dots,A_m$, $\min(A_i)=0$ and $\max(A_i)=a_i\geq0$. By the remark made at the end of Section \ref{notation_definitions} the goal is to prove the inequality
\begin{equation*}
q\left|\sum_{i=1}^{m}A_i\right|\geq\sum_{j=1}^{s}\left|\sum_{i\in S_j}A_i\right|,
\end{equation*}
where $q$ is some positive integer, the sets $S_j$ are an enumeration of the sets $S\in\mathcal{G}$ possibly with some repetition, and each $i\in[m]$ belongs to exactly $q$ of the sets $S_j$.

To start off suppose that the sets $A_1^{\prime},\dots,A_m^{\prime}$ are the sets $A_1,\dots,A_m$ translated to the left by some distances $\alpha_i\leq a_i$ such that $0\in A_i^{\prime}$, $a_i^{\prime}\geq0$ (this is immediately true since $\alpha_i\leq a_i$), and for each $j\in[s]$
\begin{equation}\label{left_volume_zero}
\left|\sum_{i\in S_j}A_i^{\prime}\cap\left(-\infty,0\right]\right|=0.
\end{equation}
At this point it may seem redundant to define these sets $A_i^{\prime}$ because the original sets $A_i$ already satisfy the above stated properties. But, later in the proof of Theorem \ref{dimension_1_equality} we will sometimes require a shift in the sets $A_i$ to achieve new sets $A_i^{\prime}$ and we must verify that certain technicalities still hold after shifting. What follows is an outline of the proof of Theorem \ref{Fractional_Superadditive}. Some of the steps are technically difficult and can be seen in full detail in \cite{Barthe1}.

\textbf{Outline of the proof of Theorem \ref{Fractional_Superadditive}:}
Since each $i\in[m]$ belongs to exactly $q$ of the sets $S_j$, there are indices 
\begin{equation*}
1\leq h_1(i)<h_2(i)<\dots<h_q(i)\leq s
\end{equation*}
such that the index $i$ belongs to only the sets $S_{h_k(i)}$. With this notation we observe that since $0\in A_i^{\prime}$ and $a_i^{\prime}\geq0$ for each $i\in[m]$, then for $1\leq k\leq q$ the following subset containments are true.
\begin{equation}\label{eq_subset}
\bigcup_{j=1}^{s}\left(\left(\sum_{i\in S_j}A_i^{\prime}+\sum_{i\in h_{k}^{-1}([1,j-1])\backslash S_j}a_i^{\prime}\right)\cap\left(\sum_{i\in h_k^{-1}([1,j-1])}a_i^{\prime};\sum_{i\in h_k^{-1}([1,j])}a_i^{\prime}\right]\right)\subset\sum_{i=1}^{m}A_i^{\prime}.
\end{equation}
Since for each $j\in[s]$ the subset containments $h_{k}^{-1}([1,j-1])\subseteq h_k^{-1}([1,j])$ are true, the sets in the union $(2)$ are disjoint. Then by taking the length (Lebesgue measure) of the left and right side of (\ref{eq_subset}) we get
\begin{equation}\label{volume_ineq1}
\sum_{j=1}^{s}\left|\left(\sum_{i\in S_j}A_i^{\prime}+\sum_{i\in h_{k}^{-1}([1,j-1])\backslash S_j}a_i^{\prime}\right)\cap\left(\sum_{i\in h_k^{-1}([1,j-1])}a_i^{\prime};\sum_{i\in h_k^{-1}([1,j])}a_i^{\prime}\right]\right|\leq\left|\sum_{i=1}^{m}A_i^{\prime}\right|.
\end{equation}
The inequalities (\ref{volume_ineq1}) hold for each $k\in[q]$. By using the translation invariance of Lebesgue measure and shifting left by $\sum_{i\in h_k^{-1}([1,j-1])\backslash S_j}a_i^{\prime}$ the inequalities (\ref{volume_ineq1}) become (full details are in \cite{Barthe1})
\begin{equation}\label{volume_ineq2}
\sum_{j=1}^{s}\left|\sum_{i\in S_j}A_i^{\prime}\cap\left(\sum_{i\in h_k^{-1}([1,j-1])\cap S_j}a_i^{\prime};\sum_{i\in h_k^{-1}([1,j])\cap S_j}a_i^{\prime}\right]\right|\leq\left|\sum_{i=1}^{m}A_i^{\prime}\right|.
\end{equation} 
Also, for fixed $j\in[s]$ as $k\in[q]$ varies, the intervals
\begin{equation}
\left(\sum_{i\in h_k^{-1}([1,j-1])\cap S_j}a_i^{\prime};\sum_{i\in h_k^{-1}([1,j])\cap S_j}a_i^{\prime}\right]
\end{equation}
are disjoint and their union over $k\in[q]$ is
\begin{equation*}
\left(0;\sum_{i\in S_j}a_i^{\prime}\right].
\end{equation*}
Then in (\ref{volume_ineq2}) by summing the left and right over $k\in[q]$ we get 
\begin{equation*}
\sum_{j=1}^{s}\left|\sum_{i\in S_j}A_i^{\prime}\cap\left(0;\sum_{i\in S_j}a_i^{\prime}\right]\right|\leq q\left|\sum_{i=1}^{m}A_i^{\prime}\right|.
\end{equation*}
Which by the assumption (\ref{left_volume_zero}) is
\begin{equation}\label{shifted_ineq}
\sum_{j=1}^{s}\left|\sum_{i\in S_j}A_i^{\prime}\right|\leq q\left|\sum_{i=1}^{m}A_i^{\prime}\right|.
\end{equation}
But by translation invariance we have that $|\sum_{i=1}^{m}A_i^{\prime}|=|\sum_{i=1}^{m}A_i|$ and for each $j\in[s]$ we have that $|\sum_{i\in S_j}A_i^{\prime}|=|\sum_{i\in S_j}A_i|$. Then (\ref{shifted_ineq}) can be written as
\begin{equation}\label{unshifted_ineq}
\sum_{j=1}^{s}\left|\sum_{i\in S_j}A\right|\leq q\left|\sum_{i=1}^{m}A_i\right|.
\end{equation}
This proves Theorem \ref{Fractional_Superadditive}.

\textbf{Concluding observations:} By observing the outline above first observe that equality in (\ref{unshifted_ineq}) holds if and only if (\ref{shifted_ineq}) is equality. Also (\ref{shifted_ineq}) is equality if and only if (\ref{volume_ineq2}) is equality for each $k\in[q]$. Moreover, (\ref{volume_ineq2}) is equality for each $k\in[q]$ if and only if for each $j\in[s]$ and $k\in[q]$
\begin{equation}\label{split_eq}
\begin{split}
&\left|\sum_{i\in S_j}A_i^{\prime}\cap\left(\sum_{i\in h_k^{-1}([1,j-1])\cap S_j}a_i^{\prime};\sum_{i\in h_k^{-1}([1,j])\cap S_j}a_i^{\prime}\right]\right|\\
&=\left|\sum_{i=1}^{m}A_i^{\prime}\cap\left(\sum_{i\in h_k^{-1}([1,j-1])\cap S_j}a_i^{\prime};\sum_{i\in h_k^{-1}([1,j])\cap S_j}a_i^{\prime}\right]\right|.
\end{split}
\end{equation}
When $j=k=1$ note that $h_{1}^{-1}(\{1\})=S_1$ and that (\ref{split_eq}) becomes
\begin{equation*}
\left|\sum_{i\in S_1}A_i^{\prime}\cap\left(0,\sum_{i\in S_1}a_i^{\prime}\right]\right|=\left|\sum_{i=1}^{m}A_i^{\prime}\cap\left(0,\sum_{i\in S_1}a_i^{\prime}\right]\right|.
\end{equation*}
Then if we show that
\begin{equation*}
\left|\sum_{i\in S_1}A_i^{\prime}\cap\left(0,\sum_{i\in S_1}a_i^{\prime}\right]\right|<\left|\sum_{i=1}^{m}A_i^{\prime}\cap\left(0,\sum_{i\in S_1}a_i^{\prime}\right]\right|
\end{equation*} 
we have shown that (\ref{split_eq}) is a strict inequality when $j=k=1$ and therefore that the inequality (\ref{unshifted_ineq}) is a strict inequality. This observation will be used frequently throughout the article.

\textbf{A remark about dimension $n\geq2$:} To prove Theorem \ref{Fractional_Superadditive} in dimension $n\geq2$ first consider the projection $\pi:\mathbb{R}^n\rightarrow\mathbb{R}$ defined by $\pi(x)=x_1$ (i.e the projection onto the first coordinate). Using translation invariance, we can assume that for the compact sets $A_i\subseteq\mathbb{R}^n$ the minimum of the projection is achieved at the origin, and that the maximum is non-negative. i.e. We assume $\pi(A_i)\subseteq [0,(a_i)_1]$, where $a_i\in A_i$ is an element which satisfies $\pi(a_i)=\max\{\pi(x):x\in A_i\}$. To prove the Theorem \ref{Fractional_Superadditive} we replace the intervals
\begin{equation*}
    \left(\sum_{i\in h_k^{-1}([1,j-1])}a_i^{\prime};\sum_{i\in h_k^{-1}([1,j])}a_i^{\prime}\right]
\end{equation*}
which appear in equation $(2)$ with the slabs
\begin{equation*}
    \pi^{-1}\left(\left(\sum_{i\in h_k^{-1}([1,j-1])}\pi(a_i^{\prime});\sum_{i\in h_k^{-1}([1,j])}\pi(a_i^{\prime})\right]\right).
\end{equation*}
The rest of the proof is the same as outlined above for the one-dimensional case. Therefore, as in the one-dimensional case, if we wish to prove that the inequality is strict, we must show that
\begin{equation}\label{slab_ineq_contradiction}
    \left|\sum_{i\in S_1}A_i^{\prime}\cap\pi^{-1}\left(\left(0,\sum_{i\in S_1}\pi(a_i^{\prime})\right]\right)\right|<\left|\sum_{i=1}^{m}A_i^{\prime}\cap\pi^{-1}\left(\left(0,\sum_{i\in S_1}\pi(a_i^{\prime})\right]\right)\right|.
\end{equation}
To summarize, in the $n$-dimensional case, we may modify the sets $A_i$ as follows. We may translate the sets $A_i$ to achieve new sets $A_i^{\prime}=A_i-\alpha_i$ as long as $0\in A_i^{\prime}$ and that for each $j$,
\begin{equation*}
    \left|\sum_{i\in S_j}A_i^{\prime}\cap\pi^{-1}\left((-\infty,0]\right)\right|=0.
\end{equation*}
We may also rotate the sets $A_i$ (but the same rotation must be done to all of them). This remark will become useful when proving Theorem \ref{dimension_n_equality}. To avoid unnecessarily tedious notation later in the paper, we will define
\begin{equation*}
    \mathcal{S}_j:=\pi^{-1}\left(\left(0,\sum_{i\in S_j}\pi(a_i^{\prime})\right]\right)
\end{equation*}
to denote the slab corresponding to the set $S_j$, and similarly $\mathcal{Z}_j$, $\mathcal{P}_j$ to denote the slabs corresponding to index sets $Z_j$ and $P_j$ respectively.

\section{Proof of Theorem \ref{dimension_1_equality}}

We will start off by observing the equality conditions when each of the sets $\sum_{i\in S}A_i$ has positive measure. This particular case will be established in Lemma \ref{positive measure}. We recall the assumptions made in the introduction that the weights in the fractional partitions are non-zero rational numbers and that a non-trivial fractional partition of $[m]$ does not contain the set $[m]$.

\begin{lemma}\label{L1}
Let $A_1,\dots,A_m$ be non-empty compact sets in $\mathbb{R}$ and let $(\mathcal{G},\beta)$ be a fractional partition of $[m]$. Then
\begin{equation*}
\left|\conv\left(\sum_{i=1}^{m}A_i\right)\right|=\sum_{S\in\mathcal{G}}\beta(S)\left|\conv\left(\sum_{i\in S}A_i\right)\right|.
\end{equation*}
\end{lemma}

\begin{proof}
To prove this lemma we will make use of the fact that for sets $A$, $B$ in $\mathbb{R}$, $\conv(A+B)=\conv(A)+\conv(B)$. Observe that 
\begin{equation*}
\begin{split}
\left|\text{conv}\left(\sum_{i=1}^{m}A_i\right)\right|&=\left|\sum_{i=1}^{m}\text{conv}(A_i)\right|\\
&=\left|\sum_{i=1}^{m}\left[\sum_{S\in\mathcal{G}:i\in S}\beta(S)\cdot\text{conv}(A_i)\right]\right|\\
&=\left|\sum_{i=1}^{m}\sum_{S\in\mathcal{G}}\beta(S)\textbf{I}_i(S)\text{conv}(A_i)\right|\\
&=\left|\sum_{S\in\mathcal{G}}\beta(S)\sum_{i\in S}\conv(A_i)\right|\\
&=\sum_{S\in\mathcal{G}}\beta(S)\left|\conv\left(\sum_{i\in S}A_i\right)\right|.
\end{split}
\end{equation*}
This completes the proof of the lemma.
\end{proof}

\begin{lemma}\label {L2}
Let $A_1,\dots,A_m$ be non-empty compact sets in $\mathbb{R}$ and let $(\mathcal{G},\beta)$ be a fractional partition over $[m]$. If each of the sets $\sum_{i\in S}A_i$ is an interval, then the set $\sum_{i=1}^{m}A_i$ is also an interval and equality is true:
\begin{equation*}
\left|\sum_{i=1}^{m}A_i\right|=\sum_{ S\in\mathcal{G}}\beta(S)\left|\sum_{i\in S}A_i\right|.
\end{equation*}
\end{lemma}

\begin{proof}
    If each $\sum_{i\in S}A_i$ is an interval, then by Lemma \ref{L1}
\begin{equation*}
\begin{split}
\left|\sum_{i=1}^{m}A_i\right|&\leq\left|\text{conv}\left(\sum_{i=1}^{m}A_i\right)\right|\\
&=\sum_{S\in\mathcal{G}}\beta(S)\left|\text{conv}\left(\sum_{i\in S}A_i\right)\right|\\
&=\sum_{S\in\mathcal{G}}\beta(S)\left|\sum_{i\in S}A_i\right|\\
&\leq\left|\sum_{i=1}^{m}A_i\right|.
\end{split}
\end{equation*}
It follows that all of the $\leq$ must be $=$ which proves that equality is true. Also we see that 
\begin{equation*}
\left|\sum_{i=1}^{m}A_i\right|=\left|\text{conv}\left(\sum_{i=1}^{m}A_i\right)\right|,
\end{equation*}
which by compactness implies that $\sum_{i=1}^{m}A_i=\conv\left(\sum_{i=1}^{m}A_i\right)$. This proves that the set $\sum_{i=1}^{m}A_i$ must be an interval.
\end{proof}

\begin{remark} In Lemma \ref{L2} the method of showing that $\sum_{i=1}^{m}A_i$ must be an interval works only because of the fact that for compact intervals $I_1, I_2\subset\mathbb{R}$ it is true that $|I_1+I_2|=|I_1|+|I_2|$. In higher dimension if $A,B$ are convex and compact, then it is not always true that $|A+B|=|A|+|B|$. This also is the reason that Lemma \ref{L1} fails and hints that maybe in dimension $n\geq2$ the equality conditions for Theorem \ref{Fractional_Superadditive} cannot have 
 a nice convexity condition like in the one-dimensional case. To prove the analogue of Lemma \ref{L2} in dimension $n\geq2$, that $\sum_{i=1}^{m}A_i$ is convex if the sets $\sum_{i\in S}A_i$ are convex for each $S\in\mathcal{G}$, there is a useful inequality involving the Schneider non-convexity index which can be used. See Section \ref{schneider} for more details.
\end{remark}

\begin{lemma}\label{translated partition}
Let $m\geq 3$ be an integer and let $(\mathcal{G},\beta)$ be a fractional partition on $[m]$ that is not the trivial partition. For $k\in(1,m)$ define a set 
\begin{equation*}
\mathcal{G^*}:=\{S\cap[k]\notin\{\varnothing,[k]\}:S\in\mathcal{G}\}.
\end{equation*}
Define the number $\gamma$ by
\begin{equation*}
\gamma:=\sum_{S\in\mathcal{G}:S\cap[k]=[k]}\beta(S).
\end{equation*}
Also assume that $\gamma<1$. Now define a function $\beta^{*}:\mathcal{G}^{*}\rightarrow[0,1]$ by
\begin{equation*}
\beta^*(T):=\frac{1}{1-\gamma}\sum_{S\in\mathcal{G}:S\cap[k]=T}\beta(S).
\end{equation*}
Then the pair $(\mathcal{G}^*,\beta^*)$ is a non-trivial fractional partition of $[k]$.
\end{lemma}

\begin{proof} 
Let $i\in[k]$. Then
\begin{equation*}
\begin{split}
\sum_{T\in\mathcal{G}^{*}:i\in T}\beta^{*}(T) &= \frac{1}{1-\gamma}\sum_{T\in\mathcal{G}^{*}:i\in T}\left[\sum_{S\in\mathcal{G}:S\cap[k]=T}\beta(S)\right]\\
&=\frac{1}{1-\gamma}\sum_{T\in\mathcal{G}^{*}}\left[\sum_{S\in\mathcal{G}:S\cap[k]=T}\beta(S)\textbf{I}_i(T)\right]\\
&=\frac{1}{1-\gamma}\sum_{T\in\mathcal{G}^{*}}\left[\sum_{S\in\mathcal{G}:S\cap[k]=T}\beta(S)\textbf{I}_i(S)\right]\\
&=\frac{1}{1-\gamma}\sum_{S\in\mathcal{G}:S\cap[k]\notin\{\varnothing,[k]\}}\beta(S)\textbf{I}_i(S)\\
&=\frac{1}{1-\gamma}\left(1-\sum_{S\in\mathcal{G}:S\cap[k]=[k]}\beta(S)\textbf{I}_i(S)\right)=1.
\end{split}
\end{equation*}
This verifies that the pair $(\mathcal{G}^{*},\beta^{*})$ is a fractional partition of $[k]$. The fact that this partition is not the trivial partition follows from the fact that each set $T\in\mathcal{G}^{*}$ is a proper subset of $[k]$ and so there must be at least two sets in $\mathcal{G}^{*}$ in order to cover all the integers in $[k]$.
\end{proof}

\begin{remark}\label{R2}
Suppose that $m\geq 3$ is an integer, that $(\mathcal{G},\beta)$ is a non-trivial fractional partition of $[m]$, and that $A_1,\dots,A_m$ are non-empty compact sets in $\mathbb{R}$. The corresponding fractional superadditive inequality is given by
\begin{equation*}
\left|\sum_{i=1}^{m}A_i\right|\geq\sum_{S\in\mathcal{G}}\beta(S)\left|\sum_{i\in S}A_i\right|.
\end{equation*}
Suppose that for some $k\in(1,m)$ exactly $k$ of the sets have $\text{card}(A_i)\geq 2$ and the remaining $m-k$ sets have $\text{card}(A_i)=1$. We will suppose that the sets $A_1,\dots,A_k$ are the sets with two or more points each, and the sets $A_{k+1},\dots,A_{m}$ are the sets which are points. Using the notation of Lemma \ref{translated partition}, we can rewrite this fractional superadditive inequality as 
\begin{equation*}
\left|\sum_{i=1}^{k}A_i\right|\geq\gamma\left|\sum_{i=1}^{k}A_i\right|+\sum_{S\in\mathcal{G}:S\cap[k]\notin\{\varnothing,[k]\}}\beta(S)\left|\sum_{i\in S\cap[k]}A_i\right|.
\end{equation*}
If $\gamma=1$, then the inequality is $|\sum_{i=1}^{k}A_i|\geq|\sum_{i=1}^{k}A_i|$. This is the fractional superadditive inequality for the trivial partition on $[k]$: $\mathcal{G}^{*}=\{[k]\}$, $\beta^{*}([k])=1$. If $\gamma<1$ we find that the inequality becomes
\begin{equation*}
\left|\sum_{i=1}^{k}A_i\right|\geq\sum_{S\in\mathcal{G}:S\cap[k]\notin\{\varnothing,[k]\}}\frac{\beta(S)}{1-\gamma}\left|\sum_{i\in S\cap[k]}A_i\right|.
\end{equation*}
By the definitions of $\mathcal{G}^{*}$ and $\beta^{*}$ this is
\begin{equation*}
\left|\sum_{i=1}^{k}A_i\right|\geq\sum_{T\in\mathcal{G}^{*}}\beta^{*}(T)\left|\sum_{i\in T}A_i\right|.
\end{equation*}
To summarize, by using the translation invariance of Lebesgue measure, we have rewritten the original inequality so that each of the sets on the right hand side $\sum_{i\in T}A_i$ contain at least two points. For future reference, the new fractional partition $(\mathcal{G}^{*},\beta^{*})$ will be called the \textit{translated fractional partition} on $[k]$ sets. That is, the new fractional partition that is obtained by accounting for translation by all of the sets $A_i$ which contain only one point.
\end{remark}

\begin{lemma} \label{positive measure}
Let $m\geq3$ be an integer, let $(\mathcal{G},\beta)$ be a fractional partition on $[m]$ that is not the trivial partition $\beta([m])=1$, and let $A_1,\dots,A_m$ be non-empty compact sets in $\mathbb{R}$ such that $|\sum_{i\in S}A_i|>0$ for any $S\in\mathcal{G}$. Then 
\begin{equation*}
\left|\sum_{i=1}^{m}A_i\right|=\sum_{S\in\mathcal{G}}\beta(S)\left|\sum_{i\in S}A_i\right|
\end{equation*}
if and only if for any $S\in\mathcal{G}$, $\sum_{i\in S}A_i$ is an interval.
\end{lemma}

\begin{proof} 
Suppose for each $S\in\mathcal{G}$ that $\sum_{i\in S}A_i$ is an interval. By Lemma \ref{L2} equality must be true.

Now assume that equality is true and assume by induction that Lemma \ref{positive measure} is true for $k<m$ sets. First, consider the case for when $k$ of the sets have $\text{card}(A_i)\geq2$ and the remaining $m-k$ sets have $\text{card}(A_i)=1$. Rearranging the sets if needed, we may assume that $A_1,\dots,A_k$ have $\text{card}(A_i)\geq2$ and that $A_{k+1},\dots,A_m$ have $\text{card}(A_i)=1$. Using the notation of Lemma \ref{translated partition} we have that 
\begin{equation*}
\gamma:=\sum_{S\in\mathcal{G}:S\cap[k]=[k]}\beta(S).
\end{equation*}
We must have that $\gamma<1$. Assume for contradiction that $\gamma=1$. Since $(\mathcal{G},\beta)$ is not the trivial partition there exists an integer $t\in[k]$ and a set $S^{(t)}\in\mathcal{G}$ such that $t\notin S^{(t)}$ (If $t$ belongs to each $S\in\mathcal{G}$, then $\sum_{S\in\mathcal{G}:t\in S}\beta(S)>1$). The set $S^{(t)}\cap[k]$ is non-empty. Otherwise the set $\sum_{i\in S^{(t)}}A_i$ would contain only point sets and it would follow that $|\sum_{i\in S^{(t)}}A_i|=0$ which is a contradiction. Choose $j\in S^{(t)}\cap[k]$. Then since $S^{(t)}\cap[k]$ is a proper subset of $[k]$ we have
\begin{equation*}
\sum_{S\in\mathcal{G}:j\in S}\beta(S)\geq\beta(S^{(t)})+\gamma>1.
\end{equation*}
This is a contradiction, so we must have $\gamma<1$. Now using Lemma \ref{translated partition} and Remark \ref{R2} the equality becomes
\begin{equation*}
\left|\sum_{i=1}^{k}A_i\right|=\sum_{T\in\mathcal{G}^{*}}\beta^{*}(T)\left|\sum_{i\in T}A_i\right|,
\end{equation*}
 where the pair $(\mathcal{G}^{*},\beta^{*})$ is a non-trivial fractional partition of $[k]$. By induction each of the sets $\sum_{i\in T}A_i$ is an interval of positive measure. By Lemma \ref{L2} the set $\sum_{i=1}^{k}A_i$ is an interval of positive measure. By translation each of the sets $\sum_{i\in S}A_i$ must be intervals of positive measure. This proves the theorem for when some of the sets have $\text{card}(A_i)=1$.

Next we will consider the case where $\text{card}(A_i)\geq2$ for all $i\in[m]$, at least one of the sets $\sum_{i\in S}A_i$ is an interval and at least one of them is not an interval. Using the assumption that the fractional partition has rational weights the equality can be written as
\begin{equation*}
q\left|\sum_{i=1}^{m}A_i\right|=\sum_{j=1}^{s}\left|\sum_{i\in S_j}A_i\right|,
\end{equation*}
where $q$ is some positive integer and the sets $S_j$ are an enumeration of the sets $S\in\mathcal{G}$ possibly with some repetition, and each $i\in [m]$ belongs to exactly $q$ of the sets $S_j$. 

\begin{claim}\label{counting_sets}
Let $\{Z_j\}$ be the sets $S_j$ such that $\sum_{i\in S_j}A_i$ is not an interval, and let $\{P_j\}$ be the sets $S_j$ such that $\sum_{i\in S_j}A_i$ is an interval. Then there exist sets $Z_j$ and $P_k$ such that $Z_j\cap P_k\subsetneq P_k$. i.e. There exists $i\in P_k$ such that $i\notin Z_j$.
\end{claim}

\begin{proof} 
Suppose for some $j_0$ that $P_k \subseteq Z_{j_0}$ for all $k$ (if not, then we are done). Since the translated partition is not trivial, there exists $i_1\in[m]$ such that $i_1\notin Z_{j_0}$. Since all the $P_k$ are contained in $Z_{j_0}$, the integer $i_1$ must belong to exactly $q$ sets $Z_j$. Now, choose some $i_2\in P_1$. If $i_2\in Z_j$ for all $j$, then $i_2$ belongs to the $q$ sets $Z_j$ that contain $i_1$, and $i_2$ belongs to $P_1$. So, $i_2$ belongs to $q+1$ sets $S_j$ which is impossible. Therefore there exists some $Z_j$ such that $i_2\notin Z_j$. This proves the claim. 
\end{proof}

Using the claim let $\sum_{i\in Z_j}A_i$ be one of the sets which is not an interval and let $\sum_{i\in P_k}A_i$ be one of the sets that is an interval which satisfies $Z_j\cap P_k\subsetneq P_k$. By rearranging the sets $S_j$ we can assume that $Z_j=S_1$. If the set $(\sum_{i=1}^{m}A_i)_{\leq\sum_{i\in S_1}a_i}$ is an interval, then since $\sum_{i\in S_1}A_i$ is not an interval,
\begin{equation*}
\left|(\sum_{i=1}^{m}A_i)_{\leq\sum_{ i\in S_1}a_i}\right|=\left|\textup{conv}\left(\sum_{i\in S_1}A_i\right)\right|>\left|\sum_{i\in S_1}A_i\right|,
\end{equation*}
which is a contradiction of (\ref{split_eq}) with $k=1$ where the sets $A_i$ have not been shifted. If the set $(\sum_{i=1}^{m}A_i)_{\leq\sum_{i\in S_1}a_i}$ is not an interval, then letting $i_0\in P_k\backslash S_1$ we have
\begin{equation*}
\begin{split}
\left|(\sum_{i=1}^{m}A_i)_{\leq\sum_{i\in S_1}a_i}\right|&=\left|\left(\textup{conv}(\sum_{i\in P_k}A_i)+\sum_{i\in [m]\backslash P_k}A_i\right)_{\leq\sum_{i\in S_1}a_i}\right|\\
&\geq\left|\left(\textup{conv}(\sum_{i\in P_k}A_i)+\sum_{(i\in [m]\backslash P_k)\cap S_1}A_i\right)_{\leq\sum_{i\in S_1}a_i}\right|\\
&\geq\left|\left(\textup{conv}(A_{i_0})+\sum_{i\in S_1}A_i\right)_{\leq\sum_{i\in S_1}a_i}\right|\\
&>\left|\sum_{i\in S_1}A_i\right|.
\end{split}
\end{equation*}
The last inequality is strict because $\text{card}(A_{i_0})\geq2$ so $\text{conv}(A_{i_0})$ is an interval with positive measure which is then added to a set with open intervals in its complement. This again proves a contradiction to (\ref{split_eq}) with $k=1$, where the sets $A_i$ have not been shifted left.

Finally, we will assume that all the sets $\sum_{i\in S_j}A_i$ are not intervals. For each $j\in[s]$ define
\begin{equation*}
\varepsilon_j:=\sup\left\{\varepsilon\geq0:\left|\sum_{i\in S_j}A_i\cap[0,\varepsilon)\right|=0\right\}.
\end{equation*}
For the following argument set $T_j:=\sum_{i\in S_j}A_i$. Suppose first for each $j\in[s]$ that $\varepsilon_j>0$. Choose $j_0\in[s]$ such that $\varepsilon_{j_0}\leq\varepsilon_j$ for each $j\in[s]$. By compactness $\varepsilon_{j_0}\in T_{j_0}$. Then $\varepsilon_{j_0}=\sum_{i\in S_{j_0}}\alpha_i$, where $\alpha_i\in A_i$ for each $i\in S_{j_0}$. Choose some $i_0\in S_{j_0}$ such that $\alpha_{i_0}>0$. Now shift the set $A_{i_0}$ to the left by $\alpha_{i_0}$ to achieve a new set $A_{i_0}^{\prime}:=A_{i_0}-\alpha_{i_0}$. Clearly after shifting left $0\in A_{i_0}^{\prime}$ and $|T_j^{\prime}\cap(-\infty,0]|=0$ for each $j\in[s]$. There exists a set $T_j$ such that $i_{0}\notin S_j$. By rearranging we may assume that this set is $T_1$. First note that $T_1+A_{i_0}^{\prime}\supseteq T_1$ because $0\in A_{i_0}^{\prime}$. Also because $A_{i_0}$ has some points to the left of $0$, $|(T_1+A_{i_0}^{\prime})\cap(0,\varepsilon_1]|>0$. Then
\begin{equation*}
\begin{split}
\left|(\sum_{i=1}^{m}A_i^{\prime})_{\leq\sum_{i\in S_1}a_i}\right|&\geq\left|(T_1+A_{i_0}^{\prime})_{\leq\sum_{i\in S_1}a_i}\right|\\
&\geq|(T_1+A_{i_0}^{\prime})\cap(0,\varepsilon_1]|+|T_1|\\
&>|T_1|.
\end{split}
\end{equation*}
This is a contradiction to (\ref{split_eq}). Next suppose for at least one $j\in[s]$ that $\varepsilon_j=0$ and for at least one $j\in[s]$ that $\varepsilon_j>0$. Choose $j_1,j_2\in[s]$ such that $\varepsilon_{j_1}=0$ and $\varepsilon_{j_2}>0$. By rearranging we may assume that $j_2=1$. By choice of $j_1$ it follows that $|T_{j_1}\cap[0,\varepsilon_1)|>0$. Then it follows that 
\begin{equation*}
\left|\sum_{i=1}^{m}A_i\cap[0,\varepsilon_1)\right|>0.
\end{equation*}
Therefore
\begin{equation*}
\left|(\sum_{i=1}^{m}A_i)_{\leq\sum_{i\in S_1}a_i}\right|\geq\left|\sum_{i=1}^{m}A_i\cap[0,\varepsilon_1)\right|+|T_1|>|T_1|.
\end{equation*}
This is a contradiction to (\ref{split_eq}). Finally suppose that $\varepsilon_j=0$ for each $j\in[s]$. Then there exists $i_0\in[m]$ such that $0$ is the limit of a decreasing sequence $\{x_k\}$ in $A_{i_0}$. Also there is $T_j$ such that $i_0\notin S_j$. By rearranging we may assume that $T_j=T_1$. Now, there exists $x\in\text{conv}(T_1)\backslash\overline{\text{PM}(T_1)}$ such that $x_L$ exists. Otherwise $T_1$ would be an interval. Choose such a number $x$. Then adding $A_{i_0}$ to $T_1$ causes a small piece of $\overline{\text{PM}(T_1)}$ to shift into the interval $(x_L,x)$. That is, $|(A_{i_0}+T_1)\cap(x_L,x)|>0$. Therefore, since $0\in A_{i_0}$ we have
\begin{equation*}
\begin{split}
\left|(\sum_{i=1}^{m}A_i)_{\leq\sum_{i\in S_1}a_i}\right|&\geq|(A_{i_0}+T_1)_{\leq\sum_{i\in S_1}a_i}|\\
&\geq|(A_{i_0}+T_1)\cap(x_L,x)|+|T_1|\\
&>|T_1|.
\end{split}
\end{equation*}
Again, this is a contradiction to (\ref{split_eq}). This completes the proof.
\end{proof}

By using Lemma \ref{positive measure} and the fact stated in Lemma \ref{L4} below that the sets $\sum_{i\in S}A_i$ with zero measure must be points, we can find the equality conditions for a general fractional partition on $m$ sets.

\begin{lemma}\label{L4}
Let $A_1,\dots,A_m$ be non-empty compact sets in $\mathbb{R}$ and let $(\mathcal{G},\beta)$ be a non-trivial fractional partition of $[m]$. Suppose that $|\sum_{i=1}^{m}A_i|>0$ and that
\begin{equation*}
\left|\sum_{i=1}^{m}A_i\right|=\sum_{S\in\mathcal{G}}\beta(S)\left|\sum_{i\in S}A_i\right|.
\end{equation*}
If $\left|\sum_{i\in S}A_i\right|=0$, then $\textup{card}\left(\sum_{i\in S}A_i\right)=1$.
\end{lemma}

\begin{proof} 
    It is enough to show that if we reduce to the translated partition that none of the sets $\sum_{i\in T}A_i$ have measure zero. So, without loss of generality we will assume that $\textup{card}(A_i)\geq2$ for each $i\in[m]$. By the rational weights assumption we have
    \begin{equation*}
        q\left|\sum_{i=1}^{m}A_i\right|=\sum_{j=1}^{s}\left|\sum_{i\in S_j}A_i\right|.
    \end{equation*}
    Let $\{P_j\}$ denote the sets $S_j$ for which $\left|\sum_{i\in S_j}A_i\right|>0$, and let $\{Z_j\}$ denote the sets $S_j$ for which $\left|\sum_{i\in S_j}A_i\right|=0$. Assume for contradiction that $\{Z_j\}\neq\varnothing$. First, suppose that one of the sets $P_j$ satisfies
    \begin{equation*}
        \left|\sum_{i\in P_j}A_i\cap(0,\varepsilon]\right|>0
    \end{equation*}
    for all $\varepsilon>0$. Assume the sets $\{S_j\}$ are ordered such that $S_1=Z_1$. Since equality is true we have (see Section $3$ equation $(8)$ with $j=k=1$)
    \begin{equation*}
        \left|\sum_{i=1}^{m}A_i\cap\left(0,\sum_{i\in S_1}a_i\right]\right|=\left|\sum_{i\in S_1}A_i\cap\left(0,\sum_{i\in S_1}a_i\right]\right|.
    \end{equation*}
    But this is a contradiction because
    \begin{equation*}
        \left|\sum_{i=1}^{m}A_i\cap\left(0,\sum_{i\in S_1}a_i\right]\right|\geq\left|\sum_{i\in P_j}A_i\cap\left(0,\sum_{i\in S_1}a_i\right]\right|>0.
    \end{equation*}
    Above we used the fact that $\sum_{i\in S_1}a_i>0$ which must be true since $\textup{card}(A_i)\geq2$ for each $i\in[m]$. Finally, assume that for each $P_j$ there exists $\varepsilon>0$ such that 
    \begin{equation*}
        \left|\sum_{i\in P_j}A_i\cap(0,\varepsilon]\right|=0.
    \end{equation*}
    For each $j$ set 
    \begin{equation*}
        \varepsilon_j:=\sup\left\{\varepsilon\geq0:\left|\sum_{i\in P_j}A_i\cap(0,\varepsilon]\right|=0\right\}.
    \end{equation*}
    By assumption, $\varepsilon_j>0$ for each $j$. In fact, all the values $\varepsilon_j$ must be the same. Otherwise, there exist $j\neq k$ such that $\varepsilon_j<\varepsilon_k$. Then
    \begin{equation*}
        \begin{split}
            \left|\sum_{i=1}^{m}A_i\cap\left(0,\sum_{i\in P_k}a_i\right]\right|&\geq\left|\left(\sum_{i=1}^{m}A_i\cap(0,\varepsilon_k]\right)\cup\left(\sum_{i=1}^{m}A_i\cap\left(\varepsilon_k,\sum_{i\in P_k}a_i\right]\right)\right|\\
            &\geq\left|\sum_{i=1}^{m}A_i\cap(0,\varepsilon_k]\right|+\left|\sum_{i\in P_k}A_i\cap\left(0,\sum_{i\in P_k}a_i\right]\right|\\
            &>\left|\sum_{i\in P_k}A_i\cap\left(0,\sum_{i\in P_k}a_i\right]\right|.
        \end{split}
    \end{equation*}
    This brings a contradiction if we set $S_1=P_k$. Now, there must exist a set $Z_j$ and a set $P_k$ such that $Z_j\cap P_k\subsetneq Z_j$. This was shown in Claim \ref{counting_sets} in the proof of Lemma \ref{positive measure}, except with the roles of $Z_j$ and $P_k$ reversed, and the sets representing the properties interval or not interval. But, with careful observation the same proof applied in Claim \ref{counting_sets} applies here. Now, choose such sets $Z_j$ and $P_k$ as just described. Then by compactness $\varepsilon_k=\sum_{i\in P_k}\alpha_i$, where $\alpha_i\in A_i$ for each $i\in P_k$. Shift all the sets $A_i$ left by $\alpha_i$ (where $\alpha_i:=0$ if $i\notin P_k$) to achieve new sets $A_i^{\prime}:=A_i-\alpha_i$. Then $0\in A_i^{\prime}$, and as noted above all the numbers $\varepsilon_j$ are the same so that $\left|\sum_{i\in S_j}A_i^{\prime}\cap(-\infty,0]\right|=0$. We also have the proper subset containment so that $\sum_{i\in Z_j}a_i^{\prime}>0$. Since we shifted left, the new value of the $\varepsilon_k$ after shifting is $\varepsilon_k^{\prime}=0$. Then
    \begin{equation*}
        \left|\sum_{i=1}^{m}A_i^{\prime}\cap\left(0,\sum_{i\in Z_j}a_i^{\prime}\right]\right|\geq\left|\sum_{i\in P_k}A_i^{\prime}\cap\left(0,\sum_{i\in Z_j}a_i^{\prime}\right]\right|>0=\left|\sum_{i\in Z_j}A_i^{\prime}\cap\left(0,\sum_{i\in Z_j}a_i^{\prime}\right]\right|.
    \end{equation*}
    If we set $Z_j=S_1$ we have a contradiction, which completes the proof of the lemma.
\end{proof}

Now we can prove the equality conditions for the one-dimensional case.

\begin{proof}[Proof of Theorem \ref{dimension_1_equality}]
Suppose that one of the three conditions holds. We will show that equality holds. If $\left|\sum_{i=1}^{m}A_i\right|=0$ or if the translated partition is trivial, then equality follows immediately. If $\left|\sum_{i=1}^{m}A_i\right|>0$, the translated partition is non-trivial, and each set $\sum_{i\in S}A_i$ is an interval, then equality follows from Lemma \ref{L2}. To prove the other direction suppose that equality holds. Suppose that $\left|\sum_{i=1}^{m}A_i\right|>0$ and that the translated partition is non-trivial. Without loss of generality we may assume that for some $k$, $\textup{card}(A_i)\geq2$ for $i=1\dots,k$, and $\textup{card}(A_i)=1$ for $i=k+1,\dots,m$. By the equality assumption we have
\begin{equation*}
    \left|\sum_{i=1}^{k}A_i\right|=\sum_{T\in\mathcal{G}^{*}}\beta^{*}(T)\left|\sum_{i\in T}A_i\right|,
\end{equation*}
where $(\mathcal{G}^*,\beta^*)$ is the translated partition on $[k]$. Since for each $T\in\mathcal{G}^*$ we have $\textup{card}\left(\sum_{i\in T}A_i\right)\geq2$, it follows from Lemma \ref{L4} that $\left|\sum_{i\in T}A_i\right|>0$ for each $T\in\mathcal{G}^*$. Then it follows from Lemma \ref{positive measure} that each set $\sum_{i\in T}A_i$ is an interval. Now, since each set $A_{k+1},\dots,A_m$ is a point, it follows that for each $S\in\mathcal{G}$ the set $\sum_{i\in S}A_i$ is an interval. This completes the proof.
\end{proof}

\section{Proof of Theorem \ref{dimension_n_equality}}

We will start off by proving Theorem \ref{dimension_n_equality} for $m=2$ sets. The reason we prove this special case separately is because the proof is enlightening as to why we might expect Theorem \ref{dimension_n_equality} to be true.

\begin{proposition} \label{true_partition}
    Let $n\geq2$ be an integer, and let $A$ and $B$ be nonempty, compact sets in $\mathbb{R}^n$. Then
    \begin{equation*}
        |A+B|=|A|+|B|
    \end{equation*}
    if and only if 
    \begin{enumerate}
        \item At least one of $A$ or $B$ is a point, or
        \item Neither $A$ nor $B$ is a point, and $|A+B|=0$.
    \end{enumerate}
\end{proposition}

\begin{proof}
    If $A$ or $B$ is a point, or if $|A+B|=0$, then equality is immediate. For the other direction suppose equality holds and that neither $A$ nor $B$ is a point. Assume for contradiction that $|A+B|>0$. By Theorem \ref{BML_inequality} and the Binomial theorem we have
    \begin{equation*}
        |A+B|\geq(|A|^{\frac{1}{n}}+|B|^{\frac{1}{n}})^n\geq|A|+|B|.
    \end{equation*}
    Note that the second inequality is actually a strict inequality of both $|A|>0$ and $|B|>0$. By the equality assumption all of the $\leq$ must be $=$. But then 
    \begin{equation*}
        |A+B|^{\frac{1}{n}}=|A|^{\frac{1}{n}}+|B|^{\frac{1}{n}}.
    \end{equation*}
    Since $A$ and $B$ both contain two or more points and $|A+B|>0$ we use the equality conditions given in Theorem \ref{BML_inequality} to deduce that $|A|>0$ and $|B|>0$.
    But then $|A+B|>|A|+|B|$ which is a contradiction. Therefore we must have $|A+B|=0$ which completes the proof.
\end{proof}

\begin{remark}
    Suppose that $(\mathcal{G},\beta)$ is a fractional partition of $[m]$ such that $\beta(S)\in\{0,1\}$ for each $S\in\mathcal{G}$. We note here that the argument of Proposition \ref{true_partition} can be extended inductively to prove that
    \begin{equation*}
        \left|\sum_{i=1}^{m}A_i\right|=\sum_{S\in\mathcal{G}}\beta(S)\left|\sum_{i\in S}A_i\right|
    \end{equation*}
    if and only if one of the two conditions is true:
    \begin{enumerate}
        \item The translated partition is trivial.
        \item The translated partition is nontrivial and $\left|\sum_{i=1}^{m}A_i\right|=0$.
    \end{enumerate}
    So, for these types of partitions the equality conditions are almost an immediate consequence of Theorem \ref{BML_inequality}.
\end{remark}

We will now prove the equality conditions for a general partition when the dimension is $n\geq2$. It should be noted that the proof is somewhat similar to the proof in one dimension, and especially makes use of techniques in the proof of Lemma \ref{L4}. The reader may wish to reference the end of Section \ref{initial_observations} for notation (especially for the definition of the slabs $\mathcal{S}_j$, $\mathcal{Z}_j$, and $\mathcal{P}_j$) and for help in understanding the approach to this proof. 

\begin{proof}[Proof of Theorem \ref{dimension_n_equality}]
    If the translated partition is trivial or if $\left|\sum_{i=1}^{m}A_i\right|=0$, then equality is immediate. For the other direction suppose that equality holds and that the translated partition is nontrivial. It is enough to prove the result on the translated partition so we will assume without loss of generality that $\textup{card}(A_i)\geq2$ for each $i\in [m]$. Assume for contradiction that $|\sum_{i=1}^{m}A_i|>0$. By the equality assumption we have
    \begin{equation*}
        q\left|\sum_{i=1}^{m}A_i\right|=\sum_{j=1}^{s}\left|\sum_{i\in S_j}A_i\right|.
    \end{equation*}
    We will prove that for some $j\in[s]$
    \begin{equation*}
        \left|\sum_{i=1}^{m}A_i\cap \mathcal{S}_j\right|>\left|\sum_{i\in S_j}A_i\cap \mathcal{S}_j\right|,
    \end{equation*}
    which by the observation made in Section \ref{initial_observations} will prove a contradiction to the equality assumption. First we will show that $|\sum_{i\in S_j}A_i|>0$ for each $j\in[s]$. If not, then let $\{Z_j\}$ denote the sets $S_j$ for which $|\sum_{i\in S_j}A_i|=0$, and let $\{P_j\}$ denote the sets $S_j$ for which $|\sum_{i\in S_j}A_i|>0$. By the same reasoning used in Claim \ref{counting_sets} in the proof of Lemma \ref{positive measure} there must exist sets $Z_j$ and $P_k$ such that $Z_j\cap P_k\subsetneq Z_j$. If the set $\sum_{i\in Z_j}A_i$ has dimension less than $n$ (i.e can be contained in a proper affine subspace of $\mathbb{R}^n$), then with the appropriate rotation we may assume that the set $\sum_{i\in Z_j}A_i$ is contained in $\pi^{-1}(0)$. Choose some $i_0\in Z_j$ such that $i_0\notin P_k$. Recall that by assumption $\textup{card}(A_{i_0})\geq2$. Then the set $A_{i_0}$ contains $0$ and also some $x\neq0$. For $\gamma\in\mathbb{R}$ define $H^+(\gamma)$ to be the translate (by $\gamma x$) of the half-space containing $x$ whose hyperplane boundary is orthogonal to $x$. Formally, we define
    \begin{equation*}
        H^+(\gamma):=\{y+\gamma x\in \mathbb{R}^n: y\cdot x\geq0\}.
    \end{equation*}
    There exists $\gamma_0\in\mathbb{R}$ such that 
    \begin{equation*}
        H^+(\gamma_0)\cap \overline{\textup{PM}}\left(\sum_{i\in P_k}A_i\right)\neq\varnothing,
    \end{equation*}
    and such that this non-empty intersection must occur entirely on the boundary of $H^+(\gamma_0)$. Therefore any shift of the set $\overline{\textup{PM}}\left(\sum_{i\in P_k}A_i\right)$ in the direction of $x$ will cause the translated set to intersect with the interior of $H^+(\gamma_0)$. Then adding $A_{i_0}$ to the set $\sum_{i\in P_k}A_i$ will cause one of the points surrounded by positive measure to move into the interior of $H^+(\gamma_0)$ and still be contained in the slab $\mathcal{P}_k$. That is,
    \begin{equation*}
        \begin{split}
            \left|\sum_{i=1}^{m}A_i\cap\mathcal{P}_k\right|&\geq\left|\left(\sum_{i\in P_k}A_i+A_{i_0}\right)\cap\mathcal{P}_k\right|\\
            &\geq\left|\sum_{i\in P_k}A_i\cap \mathcal{P}_k\right|+\left|\left(\overline{\textup{PM}}\left(\sum_{i\in P_k}A_i\right)+A_{i_0}\right)\cap H^+(\gamma_0) \cap \mathcal{P}_k\right|\\
            &>\left|\sum_{i\in P_k}A_i\cap\mathcal{P}_k\right|,
        \end{split}
    \end{equation*}
    which proves a contradiction. So, the set $\sum_{i\in Z_j}A_i$ is full dimensional (i.e cannot be contained in a proper affine subspace). Now the slab $\mathcal{Z}_j$ has positive width ($=$ the distance between the hyperplane faces). Then we can repeat the exact same argument used in Lemma \ref{L4} to find a contradiction to the equality assumption. We also note here that to be sure that the argument of Lemma \ref{L4} works properly we must make sure that the slab $\mathcal{Z}_j^{\prime}$ still has positive width after shifting. But, this must happen since the set $Z_j\backslash P_k\neq\varnothing$, and since for all $i\in Z_j\backslash P_k$, the projections $\pi(A_i)$ cannot be a single point (or else we get the same contradiction shown above). Therefore we must have $\{Z_j\}=\varnothing$, which verifies that $|\sum_{i\in S_j}A_i|>0$ for each $j\in[s]$. For $\varepsilon>0$ define
    \begin{equation*}
        B_{\varepsilon}:=\{x\in\mathbb{R}^n:\pi(x)\in[0,\varepsilon]\}.
    \end{equation*}
    For each $j\in[s]$ set
    \begin{equation*}
        \varepsilon_j:=\sup\left\{\varepsilon\geq0:\left|B_{\varepsilon}\cap\sum_{i\in S_j}A_i\right|=0\right\}.
    \end{equation*}
    By the same reasoning used in Lemma \ref{L4} all of the $\varepsilon_j$ must be the same value. Denote this common value by $\varepsilon_0$. First, assume that $\varepsilon_{0}>0$. Pick some $j_0\in[s]$. Then there exists some $\alpha\in\overline{\textup{PM}}(\sum_{i\in S_{j_0}}A_i)$ on the boundary $\{x\in\mathbb{R}^n:\pi(x)=\varepsilon_0\}$. Write $\alpha=\sum_{i\in S_{j_0}}\alpha_i$, where $\alpha_i\in A_i$. Then shift the sets $A_i$ by $\alpha_i$ to achieve new sets $A_i^{\prime}:=A_i-\alpha_i$ (if $i\notin S_{j_0}$, then set $\alpha_i=0$). Since $\varepsilon_0>0$, there exists $i_0\in S_{j_0}$ such that $\pi(\alpha_{i_0})\neq0$. There is a set $S_{j_1}$ that does not contain $i_0$. Then after shifting, the new values of $\varepsilon_j$ for the shifted sets are $\varepsilon_{j_0}^{\prime}=0<\varepsilon_{j_1}^{\prime}$. Then by the same reasoning given in Lemma \ref{L4} we have a contradiction. The only choice we are left with is $\varepsilon_{0}=0$ regardless of how the sets are rotated. We will make an additional assumption that there does not exist a rotation under which one of the hyperplane faces of a slab $\mathcal{S}_j$ contains more than one point of the set $\sum_{i\in S_j}A_i$. Otherwise we can obtain a contradiction as was shown earlier when proving that the set $\sum_{i\in Z_j}A_i$ was full dimensional. Now, choose some $j_1\in [s]$. Since $0\in\overline{\textup{PM}}(\sum_{i\in S_j}A_i)$ for each $j\in[s]$ there must exist some $i_0\in S_{j_1}$ such that $A_{i_0}$ contains a sequence $\{x_k\}$ of nonzero terms which converge to $0$. There exists a set $S_{j_2}$ which does not contain $i_{0}$. Without loss of generality we may assume that $S_{j_2}=S_1$. With this assumption we can find that the faces of $\mathcal{S}_1$ intersect the set $\sum_{i\in S_1}A_i$ as follows:
    \begin{equation*}
        \begin{split}
        \{x\in\mathbb{R}^n:\pi(x)=0\}\cap\sum_{i\in S_1}A_i&=\{0\},\\
        \{x\in\mathbb{R}^n:\pi(x)=\sum_{i\in S_1}\pi(a_i)\}\cap\sum_{i\in S_1}A_i&=\{x_{0}\}.
        \end{split}
    \end{equation*}
    Let $H$ be a hyperplane which contains $0$ and $x_0$. There are two cases to consider: Either the sequence $\{x_k\}$ has infinitely many points not contained in $H$, or the sequence has at most finitely many points not contained in $H$. Note that in dimension $n\geq3$ the second of those cases can be avoided by appropriately orienting the hyperplane $H$. Suppose first that infinitely many points of $\{x_k\}$ are not contained in $H$. Then $H$ separates $\mathbb{R}^n$ into two closed half-spaces $H^{+}$ and $H^{-}$ with boundary $H$. Let $H_{\varepsilon}^{+}$ and $H_{\varepsilon}^{-}$ be hyperplanes in $H^{+}$ and $H^{-}$ respectively that are parallel to $H$ and satisfy
    \begin{equation*}
        d\left(H_{\varepsilon}^{\pm},H\right)=\varepsilon,
    \end{equation*}
    where for sets $A$ and $B$ in $\mathbb{R}^n$ we use the notation $d(A,B):=\inf_{(a,b)\in A\times B}|a-b|$. Now, the hyperplanes $H_{\varepsilon}^{\pm}$ split $\mathbb{R}^n$ into three regions: The closed half-space $R^{+}$ bounded by $H_{\varepsilon}^{+}$ not containing $H$, the closed half-space $R^{-}$ bounded by $H_{\varepsilon}^{-}$ not containing $H$, and the intersection of the closed half spaces bounded by $H_{\varepsilon}^{\pm}$ which contain $H$, which we will denote by $R$. Now, with the above definitions fix $\varepsilon>0$ so small that $|\sum_{i\in S_1}A_i\cap R^{+}|>0$ and $|\sum_{i\in S_1}A_i\cap R^{-}|>0$ (which can be done by the assumption that all supporting hyperplanes intersect $\sum_{i\in S_1}A_i$ at a single point). By the earlier assumption, the set $\overline{\textup{PM}}(\sum_{i\in S_1}A_i)$ does not intersect the faces of $\mathcal{S}_1$ inside the regions $R^{+}$ and $R^{-}$. Now, define the numbers $\eta_{+}$ and $\eta_{-}$ by
    \begin{equation*}
        \begin{split}
            \eta_{\pm}&:=\sup\left\{\eta>0:\left|\pi^{-1}\left(\left[\sum_{i\in S_1}\pi(a_i)-\eta,\sum_{i\in S_1}\pi(a_i)\right]\right)\cap R^{\pm}\cap \overline{\textup{PM}}\left(\sum_{i\in S_1}A_i\right)\right|=0\right\}.
        \end{split}
    \end{equation*}
    By what was just said above, both $\eta_{+}$ and $\eta_{-}$ are positive. Now set $r:=\frac{1}{4}\min(\eta_{+},\eta_{-})$. Now choose $x_k\in A_{i_0}$ such that $|x_k|\leq r$. If $x_k\in \textup{int}(H^{+})$, then
    \begin{equation*}
        \begin{split}
            \left|(\sum_{i\in S_1}A_i+A_{i_0})\cap \mathcal{S}_1\right|&\geq\left|\sum_{i\in S_1}A_i\cap\mathcal{S}_1\right|+\left|\left(\sum_{i\in S_1}A_i+x_k\right)\cap\pi^{-1}\left(\left[\sum_{i\in S_1}\pi(a_i)-\eta_+,\sum_{i\in S_1}\pi(a_i)\right]\right)\cap R^{+}\right|\\
            &>\left|\sum_{i\in S_1}A_i\cap \mathcal{S}_1\right|.
        \end{split}
    \end{equation*}
    This proves a contradiction. A similar computation proves the contradiction if $x_k\in \textup{int}(H^{-})$. Now consider the other case when the dimension is $n=2$, the line $H$ through $0$ and $x_0$ is unique, and only finitely many $\{x_k\}$ are not contained in $H$. But in this case, for example, there is a line $L$, which is parallel to $H$, which passes through the set $\overline{\textup{PM}}\left(\sum_{i\in S_1}A_i\right)$ in the region $R^{+}$. Let $x_0^{\prime}$ be the intersection of $L$ with the face of $\mathcal{S}_1$ which does not contain the origin. That is,
    \begin{equation*}
        \{x\in\mathbb{R}^n:\pi(x)=\sum_{i\in S_1}\pi(a_i)\}\cap L=\{x_0^{\prime}\}.
    \end{equation*}
    Define the distance along $L$ from $x_0^{\prime}$ to the set $\overline{\textup{PM}}(\sum_{i\in S_1}A_i)$ by
    \begin{equation*}
        d:=d\left(x_0^{\prime};L\cap\overline{\textup{PM}}\left(\sum_{i\in S_1}A_i\right)\right).
    \end{equation*}
    Then $d>0$. In particular by compactness there exists $x_0^{\prime\prime}\in L\cap\overline{\textup{PM}}(\sum_{i\in S_1}A_i)$ which achieves the distance. Then choose $x_k\in H$ so that $|x_k|\leq d/2$. It follows that after adding $A_{i_0}$ to the set $\overline{\textup{PM}}(\sum_{i\in S_1}A_i)$ the point $x_0^{\prime\prime}$ is moved into the complement of the set $\overline{\textup{PM}}(\sum_{i\in S_1}A_i)$, but still is in the slab $\mathcal{S}_1$. Therefore we must have
    \begin{equation*}
        \left|\sum_{i=1}^{m}A_i\cap \mathcal{S}_1\right|\geq\left|\left(\sum_{i\in S_1}A_i+A_{i_0}\right)\cap \mathcal{S}_1\right|>\left|\sum_{i\in S_1}A_i\cap \mathcal{S}_1\right|,
    \end{equation*} 
    which is a contradiction. We have seen that in every possible case there is a contradiction, so equality is impossible and we have proved the theorem.
\end{proof}

\section{Appendix}
This appendix includes the justification of the rational weights assumption, along with an interesting observation about the Schneider non-convexity index.
\subsection{The reduction to fractional partitions with rational weights}\label{rational weights}

In proving the results of this work it has been assumed that the fractional partitions have rational weights. In general this is not true. This appendix will reiterate an argument given in \cite{Barthe1} which justifies why this assumption is sufficient.

\textbf{The reduction to rational weights:} Let $m\geq3$ be an integer. The fractional partitions of $[m]$ can be associated with sequences $(\beta_S)_{\varnothing\subsetneq S\subset[m]}$ such that for each $i\in[m]$, $\sum_{S:i\in S}\beta_S=1$. We can ignore the term $\beta_{\varnothing}$ since the empty set is irrelevant in the definition of a fractional partition. The collection of all such sequences forms a convex, compact polyhedral set $\mathcal{F}_m$:
\begin{equation*}
\mathcal{F}_m:=\left\{(\beta_S)_{\varnothing\subsetneq S\subset[m]}:\forall S, 0\leq\beta_S\leq1 \text{ and }\forall i\in[m],\sum_{S:i\in S}\beta_S=1\right\}.
\end{equation*}
Let $\beta$ be an extreme point of $\mathcal{F}_m$. There exists a set $K\subset2^{[m]}$ such that $\beta_S>0$ for any $S\in K$ and $\beta_S=0$ for $S\notin K$. It will be shown that the extreme point $\beta$ is the unique point in $\mathcal{F}_m$ such that $\beta_S>0$ for all $S\in K$ and for all $i\in[m]$, $\sum_{S\in K:i\in S}\beta_S=1$. Suppose for contradiction that there is another such $\beta^{\prime}\in \mathcal{F}_m$. Then $\beta^{\prime}_S>0$ for $S\in K$ and $\beta^{\prime}_S=0$ for $S\notin K$. Choose $\lambda>0$ so small that for each $S\in K$, $\beta_S-\lambda\beta^{\prime}_S>0$. Now define a new set of weights by $\beta^{\prime\prime}_S=(\beta_S-\lambda\beta^{\prime}_S)/(1-\lambda)$ for $S\in K$ and $\beta^{\prime\prime}_S=0$ for $S\notin K$. It can be verified that $\beta^{\prime\prime}\in \mathcal{F}_m$. Then for any $S\subset 2^{[m]}$, $\beta_S=\lambda\beta^{\prime}_S+(1-\lambda)\beta^{\prime\prime}_S$. Since $\beta$ is an extreme point of $\mathcal{F}_m$ and $\lambda\in(0,1)$ it must be that $\beta=\beta^{\prime}=\beta^{\prime\prime}$ which contradicts the assumption that $\beta$ and $\beta^{\prime}$ are distinct. Therefore, $\beta$ is the unique point of $\mathcal{F}_m$ which satisfies the system of linear equations $\sum_{S\in K:i\in S}\beta_S=1$ for each $i\in [m]$. This system has rational coefficients and by uniqueness is invertible. This implies that the non-zero weights $\beta_S$ which give the solution must be rational. Then the non-zero weights $\beta_S$ can be rewritten to have the same denominator $q$. Accounting for repeats, this implies that a fractional superadditive inequality for an extreme point $\beta$ can be written as $q|A_1+\dots+A_m|\geq\sum_{j=1}^{s}|\sum_{S_j}A|$, where the $S_j$ are an enumeration of the sets $S$ possibly with repetition and $q\geq 1$ is an integer. Moreover, each $i\in[m]$ belongs to exactly $q$ of the sets $S_j$. To see this, let $i\in[m]$. Then $i$ belongs to exactly $t$ sets in $\mathcal{G}$ (not including repetition), denote them by $S^{(1)},\dots,S^{(t)}$. Since the weights are rational with denominator $q$ we write 
\begin{equation*}
\beta(S^{(j)})=\frac{n_j}{q},
\end{equation*}
where the numbers $n_j$ are positive integers such that $n_1+\dots+n_t=q$. If $i\in S^{(j)}$, the set $S^{(j)}$ is repeated $n_j$ times, $j\in[t]$. Then $i$ belongs to at least $q$ sets (including repetition), but cannot belong to any more. It follows that there are exactly $q$ sets $S_j$ which contain $i$.

\textbf{The equivalence of Theorem \ref{dimension_1_equality} for arbitrary fractional partitions:}
Let $(\mathcal{G},\beta)$ represent some arbitrary fractional partiton of $[m]$. It is a known fact that a bounded convex polytope is the convex hull of its extreme points. Then there exist numbers $\alpha_i>0$, $i\in[N]$ with $\sum_{i=1}^{N}\alpha_i=1$ and extreme points $\beta^{(i)}$ of $\mathcal{F}_m$ such that $\beta_S=\sum_{i=1}^{N}\alpha_i\beta^{(i)}_S$ for each $S\in\mathcal{G}$. Suppose now that equality holds in the corresponding fractional superadditive inequality. Then
\begin{equation*}
\begin{split}
\left|\sum_{i=1}^{m}A_i\right|&=\sum_{S\in\mathcal{G}}\beta_S\left|\sum_{i\in S}A_i\right|\\
&=\sum_{i=1}^{N}\left(\alpha_i\sum_{S\in\mathcal{G}}\beta^{(i)}_S\left|\sum_{i\in S}A_i\right|\right).
\end{split}
\end{equation*}
From this and the fact that $|\sum_{i=1}^{m}A_i|=\sum_{i=1}^{N}\alpha_i|\sum_{i=1}^{m}A_i|$, it follows that equality holds if and only if for each $i\in[k]$
\begin{equation*}
\left|\sum_{i=1}^{m}A_i\right|=\sum_{S\in\mathcal{G}}\beta^{(i)}_S\left|\sum_{i\in S}A_i\right|.
\end{equation*}
All that is left to show is that the same equality conditions of Theorem \ref{dimension_1_equality} are also the same conditions for an arbitrary fractional partition. Using the notation of Theorem \ref{dimension_1_equality}, we only need to consider the case for when the translated partition is a non-trivial partition of $[k]$. The translated partition can be broken into a convex combination of rational partitions (which also must be non-trivial partitions of $[k]$) as described above. If a weight $\beta^{*}_S>0$, then one of the weights $\beta^{(i*)}_S$ for the rational partitions is non-zero. By Theorem 4 for rational partitions, the corresponding set $\sum_{i\in S}A_i$ must be an interval or a point depending on its volume. The reverse works the same: If the translated partition is non-trivial, the sets $\sum_{Si\in }A_i$ with positive measure are intervals, the sets with zero measure are points, then break the partition into a convex combination of rational partitions and apply Theorem \ref{dimension_1_equality}.

\subsection{Applying the Schneider non-convexity index to fractional partitions}\label{schneider}

As a technical detail in proving equality it was necessary to show (Lemma \ref{L2}) that for a fractional partition $(\mathcal{G},\beta)$ of $[m]$ and compact sets $A_1,\dots,A_m$ in $\mathbb{R}$ if for each $S\in\mathcal{G}$ the set $\sum_{i\in S}A_i$ is convex, then the set $\sum_{i=1}^{m}A_i$ is convex. This was shown by using properties of compactness and the fractional superadditive inequalities from Theorem \ref{Fractional_Superadditive}. Due to the fact that in general dimension it is not true that $|\text{conv}(A+B)|=|\text{conv}(A)|+|\text{conv}(B)|$ the technique used for proving this fact in dimension $n=1$ does not work when $n>1$. The Schneider non-convexity index \cite{Schneider1} has a strong monotonicity property \cite{Fradelizi2} which can be used to prove this technical detail when the dimension is larger than 1. For a compact set $A\subset\mathbb{R}^n$ the Schneider non-convexity index of $A$ is defined by
\begin{equation*}
c(A):=\inf\{\lambda\geq0:A+\lambda\text{conv}(A)\text{ is convex}\},
\end{equation*}
and has the property that if $A\subset\mathbb{R}^n$ is compact, $c(A)=0$ if and only if $A$ is convex. In \cite{Fradelizi2} it was shown that if $A_1,\dots,A_m$ are compact sets, and $S,T\subset [m]$ then the following inequality is true.
\begin{equation*}
c\left(\sum_{i\in S\cup T}A_i\right)\leq\max\left\{c\left(\sum_{i\in S}A_i\right),c\left(\sum_{i\in T}A_i\right)\right\}.
\end{equation*}
If $(\mathcal{G},\beta)$ is a fractional partition of $[m]$, then there exist sets $S^{(1)},\dots,S^{(t)}$ in $\mathcal{G}$ such that $S^{(1)}\cup\dots\cup S^{(t)}=[m]$. The next proposition was essentially proved in \cite{Fradelizi2}, but the result was not explicitly stated as it is here so we include it to show the interesting connection to Lemma \ref{L2}.

\begin{proposition}
Let $m\geq3$ be an integer and let $(\mathcal{G},\beta)$ be a fractional partition of $[m]$. If $A_1,\dots,A_m$ are compact sets in $\mathbb{R}^n$ such that for any $S\in\mathcal{G}$ the set $\sum_{i\in S}A_i$ is convex, then the set $\sum_{i=1}^{m}A_i$ is convex.
\end{proposition}

\begin{proof} 
By the above discussion, $c(\sum_{i\in S^{(j)}}A_i)=0$ for each $j\in[t]$. Using the above mentioned inequality it is true that for each $j<t$
\begin{equation*}
c\left(\sum_{i\in (S^{(1)}\cup\dots S^{(j)})\cup S^{(j+1)}}A_i\right)\leq\max\left\{c\left(\sum_{i\in S^{(1)}\cup\dots S^{(j)}}A_i\right),c\left(\sum_{i\in S^{(j+1)}}A_i\right)\right\}.
\end{equation*}
By induction the right hand side is equal to $0$ for each $j=1,\dots,t-1$. Setting $j=t-1$ gives
\begin{equation*}
c\left(\sum_{i=1}^{m}A_i\right)=c\left(\sum_{i\in S^{(1)}\cup\dots\cup S^{(t)}}A_i\right)=0.
\end{equation*}
This proves that the set $\sum_{i=1}^{m}A_i$ is convex.
\end{proof}

\newpage
\bibliography{Meyer_Fractional_Superadditive}
\bibliographystyle{abbrv}

\end{document}